\documentclass[12pt]{article}
\newcommand{\R}{\mathcal{R}}

\newcommand{\U}{\mathcal{U}}

\newcommand{\ita}{\textit}
\newcommand{\bo}{\textbf}
\newcommand{\A}{\mathcal{A}}

\newcommand{\M}{\mathcal{M}}

\usepackage{latexsym}
\begin{document}
\begin{center}
\textbf{\Large  Linear Programming Formulation of  Long Run Average  Optimal Control Problem}
\end{center}

\vspace{.5in}

\begin{center}
VIVEK S.\ BORKAR\footnote{Department of Electrical Engineering, Indian Institute of Technology Bombay, Powai, Mumbai 400076, India (borkar.vs@gmail.com). This work was initiated while this author was visiting Department of Mathematics, Macquarie University. Research supported in part by a J.\ C.\ Bose Fellowship from the Government of India.} AND VLADIMIR GAITSGORY\footnote{Department of Mathematics, Faculty of Science, Macquarie University, Sydney, NSW 2109, Australia (vladimir.gaitsgory@mq.edu.au). Research supported in part by  the Australian Research Council Discovery Grant
DP150100618.}
\end{center}

\vspace{.5in}

\noindent \bo{\large Abstract:} We introduce and study  the infinite dimensional  linear programming problem which along with its dual allows one to characterize the optimal value of the deterministic long-run average optimal control problem in the general  case when the latter may depend on the initial conditions of the system.\\

\noindent \bo{Key words:} long-run average optimal control; linear programming; duality; infinite horizon and vanishing discount limits

\newpage

\section{Introduction and preliminaries}
The linear programming  approach to control systems is based on the fact that the
occupational measures generated by admissible controls and the corresponding solutions of
a dynamical system satisfy certain linear equations that represent the system's dynamics in an
integral form. Such \lq\lq linearization" proved to be an efficient tool for dealing with various problems
of control, and it has been explored  extensively in both deterministic and stochastic settings (see, e.g.,  \cite{Redbook}, \cite{BhaBor}, \cite{Vivek}, \cite{BorGai}, \cite{BGQ}, \cite{F-V}, \cite{Jean}, \cite{Kurtz},   \cite{Stockbridge1}   and, respectively,   \cite{Gai8}, \cite{GQ}, \cite{GQ-1}, \cite{GR}, \cite{Goreac-Serea},  \cite{Her-Her-Lasserre}, \cite{Adelman-1}, \cite{Lass-Trelat}, \cite{QS}, \cite{Rubio}, \cite{Vinter}     as well as references therein).

In the present paper, we continue this line of research by studying  the infinite dimensional (ID) linear programming (LP) problem which along with its dual allows one to characterize the optimal value of the deterministic long-run average optimal control problem\footnote{Note that infinite time horizon optimal control problems have been traditionally studied with other (not IDLP related) techniques; see, e.g., \cite{Arisawa-1}, \cite{Arisawa-2}, \cite{Arisawa-3}, \cite{Bardi}, \cite{CHL}, \cite{GruneSIAM98}, \cite{GruneJDE98},  \cite{Sorin92}, \cite{QR-2012}, \cite{Z14}, \cite{Z06a} and references therein.} in the general  case when the latter may depend on the initial conditions of the system. Note that, while the form and the properties of the IDLP problem related to the ergodic case (that is, the case when the  optimal value is independent of the initial conditions) is well understood, the linear programming formulation of the long-run average optimal control problem in the non-ergodic case has not been discussed in the literature. In fact, a justification of such LP formulation presents a significant mathematical challenge, and (to the best of our knowledge) this is the first paper  aimed at addressing this matter.

We consider the optimal control of the system
\begin{equation}\label{e-CSO}
y'(t)=f(y(t),u(t)), \ \ \ \ \ u(t)\in U, \ \ \ \ \ t\in [0,\infty)
\end{equation}
where
  $\ f(\cdot,\cdot):  \R^m\times U \to \R^m$ is continuous in $(y,u) $
and satisfies Lipschitz condition in $y$ uniformly in $u \in U$ ($U$ is assumed to be a compact metric space). The controls
 $u(\cdot) $ are  measurable functions $u(\cdot): [0,\infty)\to U $, with the set of all controls being denoted as $\U$. Given  $u(\cdot) \in \U $ and an  initial condition $y(0)=y_0 $, the solution of system (\ref{e-CSO}) obtained with this control and this initial condition is denoted as $y(t,y_0,u) $.

 Let $Y\subset \R^m$ be a compact domain, i.e., a compact set which is the closure of its interior.
We denote by $\U_T(y_0)$, $\ \U(y_0)$  the sets of controls such that
 \begin{equation}\label{e-CSO-1}
y(t,y_0,u)\in Y \ \
\end{equation}
for any $\  t\in [0,T]$, respectively,
for any $\ t\in [0,\infty)$. (The inclusion (\ref{e-CSO-1})
can be interpreted as a state constraint.)

Consider two optimal control problems
\begin{equation}\label{Cesaro}
\frac{1}{T} \inf_{u(\cdot)\in \U_T(y_0)}\int_0^T k(y(t,y_0,u),u(t))dt:=v_T(y_0):
\end{equation}
 and
 \begin{equation}\label{Abel}
 \lambda \inf_{u(\cdot)\in \U(y_0)}\int_0^{\infty}e^{-\lambda t} k(y(t,y_0,u),u(t))dt:=h^{\lambda}(y_0),
\end{equation}
where $T>0 $, $\lambda>0 $ and
 $k(y,u):  \R^m\times U \to \R^1 $ is a continuous function. The main contribution of this paper is the introduction
 of an IDLP problem such that the limits  $\ \lim_{T\rightarrow\infty}v_T(y_0)$ and $\ \lim_{\lambda\rightarrow 0+}h^{\lambda}(y_0)$ (if they exist) are bounded from above by the optimal value of this IDLP problem and are bounded from below by the optimal value of its corresponding dual (a corollary of this being the fact that the limits $\ \lim_{T\rightarrow\infty}v_T(y_0)$ and $\ \lim_{\lambda\rightarrow 0+}h^{\lambda}(y_0)$ are equal to the optimal value of the IDLP problem provided that there is no duality gap).

 An extensive literature is devoted to matters related to the existence  and equality of the limits  $\ \lim_{T\rightarrow\infty}v_T(y_0)$ and $\ \lim_{\lambda\rightarrow 0+}h^{\lambda}(y_0)$. The ergodic case when these limits are constants (that is, when they do not depend on the initial condition $y_0$)  was studied, for example, in \cite{Redbook}, \cite{Arisawa-1}, \cite{Arisawa-2}, \cite{Arisawa-3}, \cite{Bardi}, \cite{BhaBor}, \cite{BorGai},
\cite{GQ} (see also references therein). Results  for the non-ergodic case were obtained in  \cite{BQR-2015}, \cite{GruneSIAM98}, \cite{GruneJDE98}, \cite{Sorin92},  \cite{OV-2012} and \cite{QR-2012} (of particular importance for our consideration being a nice representation for  $\ \lim_{T\rightarrow\infty}v_T(y_0)$ and $\ \lim_{\lambda\rightarrow 0+}h^{\lambda}(y_0)$ established in \cite{BQR-2015}).

In the framework of the linear programming approach, it has been shown (see \cite{GQ} and \cite{GQ-1})\footnote{Extensions of these results to degenerate diffusions appear in \cite{BhaBor}; see also \cite{Redbook}.} that in the ergodic case, the limits $\ \lim_{T\rightarrow\infty}v_T(y_0)$ and $\ \lim_{\lambda\rightarrow 0+}h^{\lambda}(y_0)$ are equal to the optimal
value of the IDLP problem
\begin{equation}\label{limits-ergodic}
k^*:= \min_{\gamma\in W}\int_{Y\times U}k(y,u)\gamma(dy,du),
\end{equation}
where
\begin{equation}\label{limits-ergodic-W}
 W:= \left\{\gamma\in \mathcal{P}(Y\times U) \ : \ \int_{Y\times U}\nabla \phi(y)^Tf(u,y)\gamma(du,dy)=0 \ \ \forall \phi(\cdot)\in C^1 \right\},
\end{equation}
with $\ \mathcal{P}(Y\times U)$ standing for the space of probability measures defined on Borel subsets of $\ Y\times U $ and  $\ C^1 $ standing for the space of continuously differentiable functions.

The IDLP problem that we are introducing in this paper is obtained by narrowing  the feasible set $W$ with the help of additional constraints allowing one to capture the dependence  of the limits   $\ \lim_{T\rightarrow\infty}\inf_{y_0\in Y}v_T(y_0)$ and $\ \lim_{\lambda\rightarrow 0+}\inf_{y_0\in Y}h^{\lambda}(y_0)$ on the initial conditions. Note that our results establishing that it is this IDLP problem and its dual that characterize  the limits of the optimal values are consistent with a celebrated result of the controlled Markov chain theory establishing
that additional constraints are needed  to characterize the limit long run average optimal value   in the non-ergodic case    (see \cite{HK-1} and \cite{HK-2}). Note that this result was obtained in the context of  Markov chains with finite state/action spaces and the corresponding finite-dimensional LP problems for which there is no duality gap (such gap being certainly a possibility in the IDLP setting; see \cite{And-1} and \cite{And-2}).

 The paper is organized as follows. After introducing key notation below, section 2 establishes lower bounds for long run average control viewed as a limiting case of finite horizon or discounted infinite horizon control problems. Section 3 derives matching upper bounds under suitable hypotheses. Together they yield the desired linear program. Section 4 considers a special case, in which there is no duality gap. Section 5 gives some longer proofs, specifically,  of a duality result and another allied result used in the foregoing.

 We conclude this section with  some  notation and definitions that are used in the sequel. First of all,
 $\ \mathcal{P}(Y\times U)$,  $\ \mathcal{M}_+(Y\times U)$ and $\ \mathcal{M}(Y\times U)$ will stand for the space of probability measures, the  space of non-negative measures and the space of all finite measures  (respectively) defined on the Borel subsets of $Y\times U$. The convergence in these  spaces will always be understood in the weak$^*$ sense,  with $\gamma^k \in \mathcal{M}(Y\times U), k =1,2,... ,$ converging to $\gamma
\in \mathcal{M}(Y\times U)$ if and only if
$\
 \lim_{k\rightarrow \infty}\int_{Y\times U} \phi(y,u) \gamma^k (dy,du) \ = \
 \int_{Y\times U} \phi(y,u) \gamma (dy,du) \
$
for any continuous $\phi(y,u): Y\times U \rightarrow \R^1$.
  The set $\mathcal{P}(Y\times U)$
  will always be treated
as a compact metric space with a metric $\rho$,
 which is  consistent
 with its weak$^*$ convergence topology.
Using this metric $\rho$,
one can define the Hausdorff metric
$\rho_H$ on the set of subsets of $\mathcal{P}(Y\times U)$ as follows: $\forall \Gamma_i \subset \mathcal{P}(Y\times U) \ , \ i=1,2 \ ,$
\vspace{-0.15cm}
\begin{equation}\label{e:intro-3}
\rho_H(\Gamma_1, \Gamma_2) := \max \{\sup_{\gamma \in
\Gamma_1} \rho(\gamma,\Gamma_2), \sup_{\gamma \in \Gamma_2}
\rho(\gamma,\Gamma_1)\}, \
\end{equation}
where $\rho(\gamma, \Gamma_i) :=
\inf_{\gamma' \in \Gamma_i} \rho(\gamma,\gamma') \ .$
Note that, although, by some abuse of terminology,  we refer to
$\rho_H(\cdot,\cdot)$ as  a metric on the set of subsets of
${\cal P} (Y \times U)$, it is, in fact, a semi-metric on this set
(note that $\rho_H(\Gamma_1, \Gamma_2)=0$ implies  $\Gamma_1
= \Gamma_2$ if  $\Gamma_1$ and $\Gamma_2$ are closed, but the equality may not be true if at least one of these sets is not closed).

 Let $u(\cdot) \in \U_T(y_0)$ and  $y(t) =
y(t,y_0,u(\cdot)), \ t\in [0,T] $. A probability measure
 $\gamma_{u(\cdot),T} \in {\cal P} (Y \times U)$ is called the
 {\it occupational measure} generated by the pair $(y(\cdot),u(\cdot) )$ on the interval $[0,T]$ if, for
  any Borel set $Q \subset Y \times
 U$,
 \begin{equation}\label{e:occup-meas-def-S}
  \gamma _{u(\cdot),T} (Q) = \frac{1}{T}\int _0 ^
 T
 1_Q (y(t),u(t)) dt ,\end{equation} where $1_Q (\cdot)$ is
the indicator function of $Q$. This definition is equivalent to
the statement that the equality
\begin{equation}\label{e:occup-meas-def-eq-S}
\int_{Y\times U} q(y,u)\gamma_{u(\cdot),T} (dy,du) = \frac{1}{T} \int _0 ^
T
 q (y(t),u(t)) dt \end{equation} is valid for any for any
$q(\cdot)\in C(Y\times U)$ (the space of continuous functions defined on $Y\times U $).

 Let $u(\cdot) \in
{\cal U}(y_0)$  and
 $y(t)=y(t,y_0,u(\cdot)),
\ t\in [0,\infty)  $. A probability measure
 $\gamma^{\lambda}_{u(\cdot)} \in {\cal P} (Y \times U)$  is called the {\it discounted occupational measure} generated
 by the pair $(y(\cdot),u(\cdot) )$ if for any Borel set $Q \subset Y \times
 U$,
\begin{equation}\label{e:occup-meas-def} \gamma ^{\lambda}_{u(\cdot)} (Q) = \lambda \int _0 ^ \infty
e^{-\lambda t} 1_Q (y(t),u(t)) dt ,\end{equation}  the latter definition
being equivalent to  the equality
\begin{equation}\label{e:occup-meas-def-eq} \int_{Y\times U} q(y,u)\gamma ^{\lambda}_{u(\cdot)} (dy,du) = \lambda \int _0 ^ \infty
e^{-\lambda t} q (y(t),u(t)) dt \end{equation}  for any
$q(\cdot)\in C(Y\times U)$.

Let $\Gamma_T(y_0)$ and $\Theta^{\lambda}(y_0) $ stand for the set of attainable occupational, respectively, discounted occupational measures:
\begin{equation}\label{e:occup-meas-def-eq-1}
\Gamma_T(y_0):= \bigcup_{u(\cdot)\in\U_T(y_0)}\{\gamma _{u(\cdot),T}\}, \ \ \ \ \ \ \ \
\Theta^{\lambda}(y_0):= \bigcup_{u(\cdot)\in\U(y_0)}\{\gamma^{\lambda}_{u(\cdot)}\}.
\end{equation}
Note that, due to (\ref{e:occup-meas-def-eq-S}) and (\ref{e:occup-meas-def-eq}), problems (\ref{Cesaro}) and (\ref{Abel}) can be reformulated in terms of occupational (resp., discounted occupational) measures as follows:
\begin{equation}\label{e:occup-meas-def-eq-2}
\inf_{\gamma\in \Gamma_T(y_0) }\int_{Y\times U}k(y,u)\gamma(dy,du) := v_T(y_0)
\end{equation}
\begin{equation}\label{e:occup-meas-def-eq-3}
\inf_{\gamma\in \Theta^{\lambda}(y_0) }\int_{Y\times U}k(y,u)\gamma(dy,du) := h^{\lambda}(y_0)
\end{equation}

\section{IDLP problems and estimates of the limit optimal values from below}\label{Section-Main-1}
Consider the IDLP problem
\begin{equation}\label{limits-non-ergodic}
\inf_{(\gamma, \xi)\in \Omega(y_0)}\int_{Y\times U}k(y,u)\gamma(dy,du):= k^*(y_0),
\end{equation}
where
$$
 \Omega(y_0):= \{(\gamma, \xi)\in \mathcal{P}(U\times Y)\times \mathcal{M}_{+}(U\times Y) \ : \ \gamma\in W,
 $$
\begin{equation}\label{non-ergodic-Omega}
\int_{Y\times U}(\phi(y_0)-\phi(y))\gamma(du,dy) + \int_{Y\times U}\nabla \phi(y)^Tf(u,y)\xi(du,dy)    =0 \ \ \forall \phi(\cdot)\in C^1 \}.
\end{equation}
Consider also the IDLP problem
\begin{equation}\label{limits-non-ergodic-dual}
 \sup_{(\mu , \psi(\cdot), \eta(\cdot) )\in \mathcal{D}}\mu :=d^*(y_0)\ \ \ \ \ \ \ \ \ \ \
\end{equation}
where $\mathcal{D}$ is the set of triplets $(\mu , \psi(\cdot), \eta(\cdot) )\in \R^1\times C^1\times C^1$
that satisfy the inequalities
\begin{equation}\label{limits-non-ergodic-dual-1}
k(y,u)+ (\psi (y_0)- \psi (y)) + \nabla \eta (y)^T f(y,u)-\mu \geq 0 \ \ \ \ \ \ \forall\ (y,u)\in Y\times U,
\end{equation}
\begin{equation}\label{limits-non-ergodic-dual-2}
 \nabla \psi (y)^T f(y,u)\geq 0 \ \ \ \ \ \ \forall\ (y,u)\in Y\times U.
\end{equation}
Problem (\ref{limits-non-ergodic-dual}) is dual to (\ref{limits-non-ergodic}) (see \cite{And-1}, \cite{And-2} and Section \ref{Section-Duality-proofs} below). In particular, the following result is valid.\\

{\bf Lemma 2.1.} {\it The optimal values of the problems (\ref{limits-non-ergodic}) and (\ref{limits-non-ergodic-dual}) are related by the inequality}
\begin{equation}\label{limits-non-ergodic-dual-4}
 k^*(y_0) \geq d^*(y_0).
\end{equation}

\bigskip

{\it Proof.} Take an arbitrary $(\gamma, \xi)\in \Omega(y_0)$ and an arbitrary  $(\mu, \psi(\cdot), \eta(\cdot))\in \mathcal{D}$. By integrating (\ref{limits-non-ergodic-dual-1})  over $\gamma$ and taking into account the fact that $\gamma\in W$, we obtain
$$
\int_{Y\times U}k(y,u)\gamma(dy,du) + \int_{Y\times U}(\psi (y_0)- \psi (y))\gamma(dy,du)\geq \mu.
$$
Also, since $(\gamma, \xi)\in \Omega(y_0)$ and since (\ref{limits-non-ergodic-dual-2}) is satisfied,
$$
\int_{Y\times U}(\psi (y_0)- \psi (y))\gamma(dy,du) = - \int_{Y\times U}\nabla\psi (y)^Tf(y,u) \xi(dy,du)\leq 0.
$$
Thus
$$
\int_{Y\times U}k(y,u)\gamma(dy,du) \geq \mu.
$$
Due to the fact that $(\gamma, \xi)\in \Omega(y_0)$ and $(\mu, \psi(\cdot), \eta(\cdot))\in \mathcal{D}  $ are arbitrary, the latter implies (\ref{limits-non-ergodic-dual-4}). $\ \Box$

\bigskip

As can be readily seen, problem (\ref{limits-non-ergodic}) can be rewritten in the form
\begin{equation}\label{limits-non-ergodic-1}
 \inf_{\gamma\in W_1(y_0)}\int_{Y\times U}k(y,u)\gamma(dy,du)=k^*(y_0),
\end{equation}
where
$$
W_1(y_0):= \{\gamma \ : \ (\gamma, \xi)\in \Omega(y_0))\}  =\{\gamma\in W \ : \ \exists \ \xi\in \mathcal{M}_{+}(Y\times U)\ \ \  {\rm such\ that}  \
$$
\begin{equation}\label{limits-non-ergodic-2}
\int_{Y\times U}(\phi(y)-\phi(y_0))\gamma(du,dy) \ =\ \int_{Y\times U}\nabla \phi(y)^Tf(u,y)\xi(du,dy) \ \forall\ \phi(\cdot)\in C^1\}.
\end{equation}
Along with problem (\ref{limits-non-ergodic-1}), let us consider the problem
\begin{equation}\label{limits-non-ergodic-3}
\min_{\gamma\in W_2(y_0)}\int_{Y\times U}k(y,u)\gamma(dy,du),
\end{equation}
where
$$
W_2(y_0):=\{\gamma\in W \ : \ \exists \ \xi_l\in \mathcal{M}_{+}(Y\times U), \ l=1,2,..., \ \  \ {\rm such\ that}
$$
\begin{equation}\label{limits-non-ergodic-4}
\int_{Y\times U}(\phi(y)-\phi(y_0))\gamma(du,dy) \ =\ \lim_{l\rightarrow\infty}\int_{Y\times U}\nabla \phi(y)^Tf(u,y)\xi_l(du,dy) \ \forall\ \phi(\cdot)\in C^1\}.
\end{equation}
It can be readily verified that the set $W_2(y_0) $ is closed (and hence compact, since $W$ is compact). Also, both $W_1(y_0) $
and $W_2(y_0)$ are convex,  with
\begin{equation}\label{limits-non-ergodic-5}
cl W_1(y_0)\subset W_2(y_0)
\end{equation}
where $cl$ stands for the closure of the corresponding set.

\bigskip

{\bf Lemma 2.2.} {\it If $W_2(y_0)\neq\emptyset$, then the optimal value of  the dual problem (\ref{limits-non-ergodic-dual}) is bounded and it is equal to the optimal value of
 problem (\ref{limits-non-ergodic-3}). That is,}
\begin{equation}\label{limits-non-ergodic-3-1}
d^*(y_0) =\min_{\gamma\in W_2(y_0)}\int_{Y\times U}k(y,u)\gamma(dy,du).
\end{equation}

\bigskip

{\it Proof.} The proof of the lemma is given in Section \ref{Section-Duality-proofs}. $\ \Box $

\bigskip

{\bf Proposition 2.3.} {\it Assume that $\ \U(y_0)\neq\emptyset$ (that is, $y_0$ belongs to the viability kernel of   (\ref{e-CSO}) in $Y$; see \cite{Aub}). Then}
\begin{equation}\label{e-main-1}
\limsup_{T\rightarrow\infty}\Gamma_T(y_0)\subset W_2(y_0),\ \ \ \ \ \ \ \ \ \ \  \liminf_{T\rightarrow\infty}v_T(y_0)\geq d^*(y_0),
\end{equation}
and
\begin{equation}\label{e-main-2}
\limsup_{\lambda\rightarrow 0}\Theta^{\lambda}(y_0)\subset W_2(y_0),\ \ \ \ \  \ \ \ \ \ \ \liminf_{\lambda\rightarrow 0}h_{\lambda}(y_0)\geq d^*(y_0).
\end{equation}

\bigskip

{\it Proof.} Note that due to our assumption that $\ \U(y_0)\neq\emptyset$, the sets
$$limsup_{T\rightarrow\infty}\Gamma_T(y_0) \ \mbox{and} \ \limsup_{\lambda\rightarrow 0}\Theta^{\lambda}(y_0) $$
are not empty. Also, as can be readily verified (see, e.g., Propositions 2.2-2.4 in \cite{GQ}),
\begin{equation}\label{e-main-3-1}
\limsup_{T\rightarrow\infty}\Gamma_T(y_0)\subset W, \ \ \ \ \ \  \limsup_{\lambda\rightarrow 0}\Theta^{\lambda}(y_0)\subset W.
\end{equation}
By (\ref{e:occup-meas-def-eq-2}) and (\ref{e:occup-meas-def-eq-3}),
$$
\liminf_{T\rightarrow\infty}v_T(y_0)=\inf\left\{\int_{Y\times U}k(y,u)\gamma(dy,du)\ | \ \gamma\in \limsup_{T\rightarrow\infty}\Gamma_T(y_0)\right\}
$$
$$
\liminf_{\lambda\rightarrow 0}h_{\lambda}(y_0)= \inf\left\{ \int_{Y\times U}k(y,u)\gamma(dy,du)\ | \ \gamma\in \limsup_{\lambda\rightarrow 0}\Theta^{\lambda}(y_0) \right\}.
$$
Therefore, by (\ref{limits-non-ergodic-3-1}), the second  relationship in (\ref{e-main-1}) and the second  relationship in (\ref{e-main-2}) follow from the corresponding first  ones.

To prove the first relationship in (\ref{e-main-1}), take any  $\ \gamma\in \limsup_{T\rightarrow\infty}\Gamma_T(y_0)$. By definition, this means that there exist  sequences
$T_l\rightarrow\infty $ and $\gamma_l\in  \Gamma_{T_l}(y_0)$ such that $\gamma_l\rightarrow\gamma $ as $l\rightarrow\infty$. The fact that the measure $\gamma_l$ belongs to the set $\Gamma_{T_l}(y_0) $ means that this measure is generated by some control $u_l(\cdot)\in \U_{T_l}(y_0) $ and the corresponding solution $y_l(t)=y(t,y_0,u_l) $ of system (\ref{e-CSO}). Consequently, for any $\phi\in C^1 $,
$$
\int_{Y\times U}(\phi(y)-\phi(y_0))\gamma_l(dy,du)= \frac{1}{T_l}\int_0^{T_l}(\phi(y_l(t))-\phi(y_0))dt
$$
\vspace{-0.4cm}
\begin{equation}\label{e-main-4}
= \frac{1}{T_l}\int_0^{T_l}\left(\int_0^t\nabla \phi(y_l(s))^Tf(y_l(s),u_l(s))ds\right)dt.
\end{equation}
Define $\xi_l\in C(Y\times U)^*$ by the equation
$$
\langle \xi_l, q \rangle = \frac{1}{T_l}\int_0^{T_l}\left(\int_0^t q(y_l(s),u_l(s))ds\right)dt \ \ \ \forall q(\cdot, \cdot)\in C(Y\times U).
$$
Note that $\ \langle \xi_l, q \rangle\geq 0 $ if $q(\cdot, \cdot)\geq 0 $. Hence by Riesz theorem (Theorem 4.3.9, p.\ 181 in \cite{Ash}) there exists $\xi_l\in \M_+(Y\times U) $ such that
$$
\langle \xi_l, q \rangle = \int_{Y\times U} q(y,u)\xi_l(dy,du)\ \ \ \forall q(\cdot, \cdot)\in C(Y\times U).
$$
Taking these relationships into consideration, one can rewrite (\ref{e-main-4}) in the form
$$
\int_{Y\times U}(\phi(y)-\phi(y_0))\gamma_l(dy,du) = \int_{Y\times U} \nabla \phi(y)^Tf(y,u)\xi_l(dy,du).
$$
Passing to the limit in the expression above, one obtains
$$
\int_{Y\times U}(\phi(y)-\phi(y_0))\gamma(dy,du) = \lim_{l\rightarrow\infty}\int_{Y\times U} \nabla \phi(y)^Tf(y,u)\xi_l(dy,du).
$$
Since by (\ref{e-main-3-1}), $\gamma\in W$, the latter implies that $\gamma\in W_2(y_0) $. Thus the first relationship in (\ref{e-main-1}) is established.

To prove the first relationship in (\ref{e-main-2}), note that
\begin{equation}\label{e-main-5-0}
\Theta^{\lambda}(y_0)\subset W(\lambda, y_0),
\end{equation}
where
$$
W(\lambda, y_0)= \{\gamma\in \mathcal{P}(Y\times U) \ :
$$
\begin{equation}\label{e-main-5}
 \ \int_{U\times Y}\left(\nabla \phi(y)^Tf(u,y)+ \lambda(\phi(y_0)-\phi(y)\right)\gamma(dy,du)=0 \ \ \forall \phi(\cdot)\in C^1\};
\end{equation}
see, e.g., Proposition 2.2 in \cite{GQ}. (In fact under certain non-restrictive conditions, the closed convex hull of $\Theta^{\lambda}(y_0)$ is equal to $W(\lambda, y_0)$; see \cite{GQ} and \cite{GQ-1}.)

By (\ref{e-main-5-0}), to prove that
$\ \limsup_{\lambda\rightarrow 0}\Theta^{\lambda}(y_0)\subset W_2(y_0) $, it is sufficient to prove that
\begin{equation}\label{e-main-6}
\limsup_{\lambda\rightarrow 0}W(\lambda, y_0)\subset W_2(y_0) .
\end{equation}
Note that it can be readily verified (see, e.g., Lemma 2.4 in  \cite{GQ}) that
\begin{equation}\label{e-main-7}
\limsup_{\lambda\rightarrow 0}W(\lambda, y_0)\subset W .
\end{equation}
Take now an arbitrary $\gamma\in \limsup_{\lambda\rightarrow 0}W(\lambda, y_0) $. By definition, it means that there exists sequences
$\lambda_l\rightarrow 0 $ and $\ \gamma_l\in W(\lambda_l, y_0) $ such that $\ \gamma_l\rightarrow \gamma $ as $l\rightarrow\infty$.
Since $\gamma_l\in W(\lambda_l, y_0) $, we have
\begin{equation}\label{e-main-8}
\int_{U\times Y}(\phi(y)-\phi(y_0))\gamma_l(dy,du)= \frac{1}{\lambda_l}\int_{U\times Y}\nabla \phi(y)^Tf(u,y)\gamma_l(dy,du).
\end{equation}
Define $\xi_l\in C(Y\times U)^*$ by the equation
$$
\langle \xi_l, q \rangle = \frac{1}{\lambda_l}\int_{U\times Y}q(y,u)\gamma_l(dy,du) \ \ \ \forall q(\cdot, \cdot)\in C(Y\times U).
$$
Note that $\ \langle \xi_l, q \rangle\geq 0 $ if $q(\cdot, \cdot)\geq 0 $. Hence (see Theorem 4.3.9, p. 181 in \cite{Ash}) there exists $\xi_l\in \M_+(Y\times U) $ such that
$$
\langle \xi_l, q \rangle = \int_{Y\times U} q(y,u)\xi_l(dy,du)\ \ \ \forall q(\cdot, \cdot)\in C(Y\times U).
$$
Thus (\ref{e-main-8}) can be rewritten in the form
$$
\int_{U\times Y}(\phi(y)-\phi(y_0))\gamma_l(dy,du)= \int_{U\times Y}\nabla \phi(y)^Tf(u,y)\xi_l(dy,du).
$$
Passing to the limit in this expression, one obtains
$$
\int_{U\times Y}(\phi(y)-\phi(y_0))\gamma(dy,du)= \lim_{l\rightarrow\infty}\int_{U\times Y}\nabla \phi(y)^Tf(u,y)\xi_l(dy,du).
$$
The latter, along with the fact that $\gamma\in W $ (see (\ref{e-main-7})), establish the validity of the first relationship in (\ref{e-main-2}).  $\ \Box$

\bigskip

Let $u_{\mathcal{T}}(\cdot)\in \U(y_0)$ be $\mathcal{T}$-periodic (for some $\mathcal{T}>0$). Assume that corresponding to this periodic control, there exists a $\mathcal{T}$-periodic solution $y_{\mathcal{T}}(t)= y(t, y_0, u_{\mathcal{T}}) $ of system (\ref{e-CSO}). The pair $(y_{\mathcal{T}}(\cdot), u_{\mathcal{T}}(\cdot)) $ will be referred to as a {\it $y_0 $-admissible  $\mathcal{T}$-periodic pair}. Consider the optimal control problem (this being commonly referred to as periodic optimization problem)
\begin{equation}\label{e-main-8-1}
\inf_{\mathcal{T},\left(y_{\mathcal{T}}(\cdot),u_{\mathcal{T}}(\cdot)\right)  }\left\{\frac{1}{\mathcal{T}}\int_0^{\mathcal{T}} k(y_{\mathcal{T}}(t),u_{\mathcal{T}}(t))dt \right\} :=v_{per}(y_0),
\end{equation}
where $inf$ is over all $\mathcal{T}>0$ and over all $y_0$-admissible $\mathcal{T}$-periodic pairs $(y_{\mathcal{T}}(\cdot), u_{\mathcal{T}}(\cdot)) $. Similarly to (\ref{e:occup-meas-def-eq-2}), this problem can be reformulated in terms of occupational measures
\begin{equation}\label{e:occup-meas-def-eq-per}
\inf_{\gamma\in \Gamma_{per}(y_0) }\int_{Y\times U}k(y,u)\gamma(dy,du)=v_{per}(y_0),
\end{equation}
where $\Gamma_{per}(y_0)$ is the set of occupational measures generated by all $y_0 $-admissible  periodic pairs. Note that
\begin{equation}\label{e:occup-meas-def-eq-per-11}
\Gamma_{per}(y_0)\subset \limsup_{T\rightarrow\infty}\Gamma_T(y_0)
\end{equation}
and, therefore,
\begin{equation}\label{e:occup-meas-def-eq-per-1}
v_{per}(y_0)\geq \liminf_{T\rightarrow\infty}v_T(y_0).
\end{equation}

\bigskip

{\bf Proposition 2.4.}{\it The following relationships are valid:}
\begin{equation}\label{e-main-8-2}
\Gamma_{per}(y_0)\subset W_1(y_0), \ \ \ \ \ \ v_{per}(y_0)\geq k^*(y_0).
\end{equation}

\bigskip

{\it Proof.} Due to (\ref{limits-non-ergodic-1}) and (\ref{e:occup-meas-def-eq-per}), it is sufficient to prove only the first relationship. Note that from (\ref{e-main-3-1}) and (\ref{e:occup-meas-def-eq-per-11}) it follows that
\begin{equation}\label{e-main-8-2-5}
\Gamma_{per}(y_0)\subset W.
\end{equation}
Take now an arbitrary $\gamma\in \Gamma_{per}(y_0)$. By definition, it means that $\gamma$ is generated by some $y_0 $-admissible  $\mathcal{T}$-periodic pair $(y_{\mathcal{T}}(\cdot), u_{\mathcal{T}}(\cdot)) $. That is, for any continuous function $q(y,u) $,
$$
\int_{(y,u)\in Y\times U}q(y,u)\gamma(dy,du)= \frac{1}{\mathcal{T}}\int_0^{\mathcal{T}}q(y_{\mathcal{T}}(t), u_{\mathcal{T}}(t))dt
$$
Consequently, for any $\phi\in C^1 $,
$$
\int_{Y\times U}(\phi(y)-\phi(y_0))\gamma(dy,du)= \frac{1}{\mathcal{T}}\int_0^{\mathcal{T}}(\phi(y_{\mathcal{T}}(t))-\phi(y_0))dt
$$
\vspace{-0.4cm}
\begin{equation}\label{e-main-4-per}
= \frac{1}{\mathcal{T}}\int_0^{\mathcal{T}}\left(\int_0^t\nabla \phi(y_{\mathcal{T}}(s))^Tf(y_{\mathcal{T}}(s),u_{\mathcal{T}}(s))ds\right)dt.
\end{equation}
Define $\xi\in C^*(Y\times U)$ by the equation
$$
\langle \xi, q \rangle = \frac{1}{\mathcal{T}}\int_0^{\mathcal{T}}\left(\int_0^{\mathcal{T}} q(y_{\mathcal{T}}(s),u_{\mathcal{T}}(s))ds\right)dt \ \ \ \forall q(\cdot, \cdot)\in C(Y\times U).
$$
As one can see, $\ \langle \xi, q \rangle\geq 0 $ if $q(\cdot, \cdot)\geq 0 $. Hence, by Riesz theorem (see Theorem 4.3.9, p. 181 in \cite{Ash}) there exists $\xi\in \M_+(Y\times U) $ such that
$$
\langle \xi, q \rangle = \int_{Y\times U} q(y,u)\xi(dy,du)\ \ \ \forall q(\cdot, \cdot)\in C(Y\times U).
$$
Taking these relationships into consideration, one can rewrite (\ref{e-main-4-per}) in the form
$$
\int_{Y\times U}(\phi(y)-\phi(y_0))\gamma(dy,du) = \int_{Y\times U} \nabla \phi(y)^Tf(y,u)\xi(dy,du).
$$
Since $\gamma\in W$ (by (\ref{e-main-8-2-5})), the latter implies that $\gamma\in W_1(y_0) $. Thus the first relationship in (\ref{e-main-8-2}) is established. $\ \Box $

\bigskip

{\bf Corollary 2.5.} {\it If
\begin{equation}\label{e-main-8-3}
 v_{per}(y_0) = \liminf_{T\rightarrow\infty}v_T(y_0),
\end{equation}
then}
\begin{equation}\label{e-main-8-4}
\liminf_{T\rightarrow\infty}v_T(y_0)\geq k^*(y_0).
\end{equation}
\section{Estimate of the limit optimal values from above}\label{Section-Main-2}

Let is introduce the following  assumptions:

\bigskip

ASSUMPTION I.  The set $Y$ is forward invariant with respect to solutions of system (\ref{e-CSO}). That is,  $\U(y_0)= \U\ \forall y_0\in Y$.

\medskip

ASSUMPTION II.  For any $y_0\in \mathcal{N}$, where $\mathcal{N}$ is an open neighbourhood of $Y$, there exists a limit
\begin{equation}\label{eq-est-from-above-0}
\lim_{T\rightarrow\infty}v_T(y_0):=v^*(y_0),
\end{equation}
the convergence being uniform with respect to $y_0\in \mathcal{N}$ and the limit function $v^*(\cdot) $ being Lipschitz continuous on $\mathcal{N}$.
 (Note that from this assumption it follows  that
$\lim_{\lambda\rightarrow 0}h^{\lambda}(y_0)=v^*(y_0) $; see \cite{OV-2012}.)

\medskip

{\bf Proposition 3.1.} {\it Let problem (\ref{limits-non-ergodic}) be consistent (that is, $\Omega(y_0)\neq\emptyset$)\footnote{Note that $\Omega(y_0)\neq\emptyset$ if $\Gamma_{per}\neq\emptyset$; see Proposition 2.4. } and let   Assumptions I and II be valid. Then the limit optimal value $v^*(y_0)$ is less or equal than the optimal value of the IDLP problem (\ref{limits-non-ergodic}). That is,}
\begin{equation}\label{eq-est-from-above-3}
v^*(y_0)\leq k^*(y_0)\ \ \ \forall \ y_0\in Y.
\end{equation}

\bigskip

{\it Proof.}
In \cite{BQR-2015} it has been established that, under Assumptions I and II, the limit  function $v^*(y_0) $ allows  the following representation for $y_0\in Y $:
\begin{equation}\label{eq-est-from-above-1}
v^*(y_0)=\sup_{w(\cdot)\in \mathcal{H}}w(y_0)\ \ \ \ \ \forall \ y_0\in Y,
\end{equation}
where $\mathcal{H}$ is the set of Lipschitz continuous on $\mathcal{N}$ functions\footnote{In \cite{BQR-2015}, the functions in   $\mathcal{H}$ were assumed to be just continuous (not Lipschitz continuous). However, if $v^*(\cdot) $ is Lipschitz continuous, then the representation (\ref{eq-est-from-above-1}) is valid with $\mathcal{H}$ consisting of Lipschitz continuous functions.}  such that $w(\cdot)\in \mathcal{H} $ if and only if:

$(i)$ For any $y_0\in \mathcal{N} $  and any $u(\cdot)\in \U $, the function $w(y(t,y_0,u))$ is  nondecreasing in $t$  on any interval $t\in [0,T)\ (T>0) $ such that
$y(t,y_0,u)\in \mathcal{N}\ \forall t\in [0,T]$;

\medskip

$(ii) $ The following inequality is valid:
\begin{equation}\label{eq-est-from-above-2}
\int_{Y\times U}w(y)\gamma(du,dy)\leq \int_{Y\times U}k(y,u)\gamma(du,dy)\ \ \ \ \forall \ \gamma\in W.
\end{equation}

Since any function $w(\cdot)\in \mathcal{H}$ is  Lipschitz continuous, it is almost everywhere differentiable on $\mathcal{N}$
(by Rademacher  theorem; see \cite{EG92}). Moreover, due to the property $(i)$ of the set $\mathcal{H}$,
\begin{equation}\label{eq-est-from-above-3-1}
\nabla w(y)^T f(y,u)\geq 0 \ \ \ \forall \ u\in U
\end{equation}
at any point $y\in \mathcal{N}$ where $\nabla w(y)$ exists. Let the set  $D^*_{w}(y) $ be defined as follows
\begin{equation}\label{eq-est-from-above-3-2}
D^*_{w}(y):= \{p\in \R^m \ | \ p = \lim_{l\rightarrow\infty}\nabla w(y_l) \ \ {\rm for \ some}\ \ y_l\rightarrow y \}.
\end{equation}
The set $D^*_{w}(y) $ is non-empty and compact  for any $y\in \mathcal{N} $, and  by (\ref{eq-est-from-above-3-1}),
\begin{equation}\label{eq-est-from-above-3-3}
 p^T f(y,u)\geq 0 \ \ \ \forall \ u\in U , \ \  \forall \ p\in D^*_{w}(y) .
\end{equation}
According to a well-known result in non-smooth analysis,
$$
co D^*_{w}(y)= \partial w(y)
$$
where $\partial w(y) $ is  the generalized Clarke's gradient (see, e.g., p. 63 in \cite{Bardi}). Therefore,
\begin{equation}\label{eq-est-from-above-3-4}
 p^T f(y,u)\geq 0 \ \ \ \forall \ u\in U , \ \  \forall \ p\in \partial w(y)
\end{equation}
for any $y\in \mathcal{N} $.
Let now $\epsilon > 0$ be small enough so that
\begin{equation}\label{eq-est-from-above-3-5-1}
Y + \epsilon B\subset \mathcal{N},
\end{equation}
where $B$ is the open unit ball in $\R^m $. Due to (\ref{eq-est-from-above-1}), there exists
$w_{\epsilon}(\cdot)\in \mathcal{H}  $  such that
\begin{equation}\label{eq-est-from-above-4}
w_{\epsilon}(y_0)\geq v^*(y_0) - \epsilon.
\end{equation}
By Theorem 2.2 in \cite{CZAR}, there exists $\psi_{\epsilon}(\cdot)\in C^1 $ such that
\begin{equation}\label{eq-est-from-above-5}
\max_{y\in Y}|\psi_{\epsilon}(y) - w_{\epsilon}(y)|\leq \epsilon,
\end{equation}
and such that
\begin{equation}\label{u-2}
\displaystyle{ \nabla \psi _\epsilon (y) \in \bigcup _{y' \in y
+\epsilon B} \partial \psi (y ' ) + \epsilon B\ \ \ \ \forall \ y\in Y.
}\end{equation}
Note that from (\ref{eq-est-from-above-4}) and (\ref{eq-est-from-above-5}) it follows that
\begin{equation}\label{eq-est-from-above-6}
\psi_{\epsilon}(y_0)\geq v^*(y_0) - 2 \epsilon
\end{equation}
and also (since $w_{\epsilon}(\cdot)\in \mathcal{H}  $; see (\ref{eq-est-from-above-2})) that
\begin{equation}\label{eq-est-from-above-7}
\int_{Y\times U}\psi_{\epsilon}(y)\gamma(du,dy)\leq \int_{Y\times U}(k(y,u)+\epsilon)\gamma(du,dy)\ \ \ \ \forall \ \gamma\in W.
\end{equation}
From (\ref{u-2}), on the other hand, it follows that, for an arbitrary $y \in Y $, there
exist $\ y _{\epsilon} \in y + \epsilon B $, $ \ d_{\epsilon} \in
\partial \psi (y _{\epsilon}) $ and  $h_{\epsilon} \in \epsilon B$ such that
\begin{equation}\label{eq-est-from-above-8-1}
\nabla \psi_{\epsilon} (y)= d_{\epsilon}+ h_{\epsilon},
\end{equation}
with $\ d_{\epsilon}^Tf(y,u)\geq 0\ \forall u\in U$ (due to (\ref{eq-est-from-above-3-4}) and (\ref{eq-est-from-above-3-5-1}))
and with $\ || h_{\epsilon} ||\leq \epsilon  $. Therefore,
$$
\nabla \psi_{\epsilon}(y)^Tf(y,u)=  d_{\epsilon}^T f(y,u) + h_{\epsilon}^T f(y,u)\geq -\epsilon ||f(y,u)||.
$$
Consequently,
\begin{equation}\label{eq-est-from-above-9}
\nabla \psi_{\epsilon}(y)^Tf(y,u)\geq - \epsilon M \ \ \forall (y,u)\in Y\times U, \ \ \ {\rm where} \ \ M:=\max_{(y,u)\in Y\times U}||f(y,u)||.
\end{equation}

Let us now rewrite inequality (\ref{eq-est-from-above-7}) in the form
$$
 \int_{Y\times U}(k(y,u)+\epsilon - \psi_{\epsilon}(y))\gamma(du,dy) \geq 0\ \ \ \ \forall \ \gamma\in W
$$
implying that
\begin{equation}\label{eq-est-from-above-10}
\min_{\gamma\in W} \int_{Y\times U}(k(y,u)+\epsilon - \psi_{\epsilon}(y))\gamma(du,dy) \geq 0.
\end{equation}
The problem on the left-hand-side of (\ref{eq-est-from-above-10}) is of the IDLP class, the dual of which is
\begin{equation}\label{eq-est-from-above-11}
\sup_{\eta(\cdot)\in C^1}\min_{(y,u)\in Y\times U}\left\{  k(y,u)+\epsilon - \psi_{\epsilon}(y) + \nabla \eta(y)^T f(y,u)\right\}.
\end{equation}
The optimal values of the former and the latter are equal (see, Theorem 4.1 in \cite{FinGaiLeb} or  Theorem 3.1 in \cite{GQ}). Therefore, the
inequality (\ref{eq-est-from-above-10}) can be rewritten in the form
\begin{equation}\label{eq-est-from-above-12}
\sup_{\eta(\cdot)\in C^1}\min_{(y,u)\in Y\times U}\left\{  k(y,u)+\epsilon - \psi_{\epsilon}(y) + \nabla \eta(y)^T f(y,u)\right\} \geq 0.
\end{equation}
From (\ref{eq-est-from-above-12}), it follows that there exists a function $\eta_{\epsilon}(\cdot)\in C^1 $ such that
$$
\min_{(y,u)\in Y\times U}\left\{  k(y,u)+\epsilon - \psi_{\epsilon}(y) + \nabla \eta_{\epsilon}(y)^T f(y,u)\right\} \geq -\epsilon.
$$
That is,
\begin{equation}\label{eq-est-from-above-13}
 k(y,u) - \psi_{\epsilon}(y) + \nabla \eta_{\epsilon}(y)^T f(y,u)  \geq -2\epsilon\ \ \ \forall (y,u)\in Y\times U .
\end{equation}
Consider the following IDLP problem
\begin{equation}\label{eq-est-from-above-14}
\sup_{(\psi(\cdot),\eta(\cdot))\in \mathcal{Q}(\epsilon)}\psi(y_0):=\bar d^*(y_0,\epsilon)
\end{equation}
where $\mathcal{Q}(\epsilon)$ is the set of pairs $(\psi(\cdot), \eta(\cdot) )\in  C^1\times C^1$
that satisfy the inequalities
\begin{equation}\label{limits-non-ergodic-dual-1-eps}
k(y,u)- \psi (y) + \nabla \eta (y)^T f(y,u) \geq -2\epsilon \ \ \ \ \ \ \forall\ (y,u)\in Y\times U,
\end{equation}
\begin{equation}\label{limits-non-ergodic-dual-2-eps}
 \nabla \psi (y)^T f(y,u)\geq  -M\epsilon \ \ \ \ \ \ \forall\ (y,u)\in Y\times U.
\end{equation}
Note that, by (\ref{eq-est-from-above-9}) and (\ref{eq-est-from-above-13}), $(\psi_{\epsilon}(\cdot), \eta_{\epsilon}(\cdot))\in \mathcal{Q}(\epsilon)$.
Consequently (and also due to (\ref{eq-est-from-above-6})),
\begin{equation}\label{eq-est-from-above-15}
\bar d^*(y_0,\epsilon)\geq \psi_{\epsilon}(y_0)\geq v^*(y_0) - 2 \epsilon .
\end{equation}
Consider also the problem
\begin{equation}\label{limits-non-ergodic-dual-eps-d}
 \sup_{(\mu , \psi(\cdot), \eta(\cdot) )\in \mathcal{D}(\epsilon)}\mu :=d^*(y_0,\epsilon)\ \ \ \ \ \ \ \ \ \ \
\end{equation}
where $\mathcal{D}(\epsilon)$ is the set of triplets $(\mu , \psi(\cdot), \eta(\cdot) )\in \R^1\times C^1\times C^1$
that satisfy the inequalities
\begin{equation}\label{limits-non-ergodic-dual-eps-1-d}
k(y,u)+ (\psi (y_0)- \psi (y)) + \nabla \eta (y)^T f(y,u)-\mu \geq  - 2 \epsilon \ \ \ \ \ \ \forall\ (y,u)\in Y\times U,
\end{equation}
\begin{equation}\label{limits-non-ergodic-dual-2-eps-d}
 \nabla \psi (y)^T f(y,u)\geq  - M \epsilon \ \ \ \ \ \ \forall\ (y,u)\in Y\times U.
\end{equation}
Let us show that the optimal values of (\ref{eq-est-from-above-14}) and (\ref{limits-non-ergodic-dual-eps-d}) are equal.
That is,
\begin{equation}\label{eq-est-from-above-165}
\bar d^*(y_0,\epsilon)= d^*(y_0,\epsilon) .
\end{equation}
Firstly, note that $\bar d^*(y_0,\epsilon)\leq d^*(y_0,\epsilon) $ (since, for any pair $(\psi(\cdot), \eta(\cdot))\in \mathcal{Q}(\epsilon)  $, the triplet $(\mu, \psi(\cdot), \eta(\cdot))\in \mathcal{D}(\epsilon)  $, where $\mu=\psi(y_0)$;
see (\ref{limits-non-ergodic-dual-1-eps})-(\ref{limits-non-ergodic-dual-2-eps}) and (\ref{limits-non-ergodic-dual-eps-1-d})-(\ref{limits-non-ergodic-dual-2-eps-d})). Let us prove the converse inequality. Let
a triplet $(\mu', \psi'(\cdot), \eta'(\cdot))\in \mathcal{D}(\epsilon)  $ be  such that $\mu'\geq d^*(y_0,\epsilon)- \delta $,
with $\delta > 0 $ being arbitrarily small. Then the pair $(\tilde \psi'(\cdot), \eta'(\cdot))\in \mathcal{Q}(\epsilon) $, where
$\tilde \psi'(y):= \psi'(y)- \psi'(y_0)+\mu' $. Since $\tilde \psi'(y_0) = \mu' $, it leads to the inequality $\bar d^*(y_0,\epsilon)\geq d^*(y_0,\epsilon)- \delta$, and, consequently, to the inequality $\bar d^*(y_0,\epsilon)\geq d^*(y_0,\epsilon) $ (since $\delta > 0$ is arbitrarily small). Thus (\ref{eq-est-from-above-165}) is proved.

Problem (\ref{limits-non-ergodic-dual-eps-d}) is dual to the IDLP problem
\begin{equation}\label{limits-non-ergodic-pert-eps}
\inf_{(\gamma, \xi)\in \Omega(y_0)}\left\{\int_{Y\times U}(k(y,u)+ 2\epsilon)\gamma(dy,du)+ M\epsilon \int_{Y\times U} \xi(dy,du) \right\}:= k^*(y_0, \epsilon).
\end{equation}
As established by Lemma  5.1 (see Section \ref{Section-Duality-proofs}), the optimal values of (\ref{limits-non-ergodic-dual-eps-d}) and (\ref{limits-non-ergodic-pert-eps}) are equal for any $\epsilon > 0$. That is,
\begin{equation}\label{eq-est-from-above-16}
d^*(y_0, \epsilon) = k^*(y_0, \epsilon)\ \ \ \forall \ \epsilon > 0.
\end{equation}
Therefore, by (\ref{eq-est-from-above-15}) and  (\ref{eq-est-from-above-165}),
\begin{equation}\label{eq-est-from-above-17}
k^*(y_0, \epsilon)\geq v^*(y_0) - 2\epsilon.
\end{equation}
Note that problem (\ref{limits-non-ergodic-pert-eps}) is a perturbed version of problem (\ref{limits-non-ergodic}) and that (due to (\ref{eq-est-from-above-17})), to prove
(\ref{eq-est-from-above-3}) it is sufficient to establish that
\begin{equation}\label{eq-est-from-above-18}
\lim_{\epsilon\rightarrow 0}k^*(y_0, \epsilon)= k^*(y_0).
\end{equation}
As can be easily seen, $\ k^*(y_0, \epsilon)$ is a decreasing function of $\epsilon $, and $\ k^*(y_0, \epsilon)\geq k^*(y_0) \ \forall \ \epsilon > 0$. Hence,
$$
\lim_{\epsilon\rightarrow 0}k^*(y_0, \epsilon)\geq k^*(y_0).
$$
Let $\delta > 0 $ be arbitrary small and $(\gamma', \xi')\in \Omega(y_0) $
be $\delta$-near-optimal for (\ref{limits-non-ergodic}). That is,
$$
 \int_{Y\times U}k(y,u)\gamma'(dy,du)\leq k^*(y_0) + \delta.
$$
Then
$$
k^*(y_0, \epsilon)\leq \int_{Y\times U}(k(y,u)+ 2\epsilon)\gamma'(dy,du)+ M\epsilon \int_{Y\times U} \xi'(dy,du)
$$
$$
\leq  k^*(y_0) + \delta +2\epsilon \int_{Y\times U}\gamma'(dy,du) + M\epsilon \int_{Y\times U}\xi'(dy,du),
$$
$$
\Rightarrow\ \ \ \ \lim_{\epsilon\rightarrow 0}k^*(y_0, \epsilon)\leq k^*(y_0) + \delta \ \ \ \Rightarrow \ \ \ \lim_{\epsilon\rightarrow 0}k^*(y_0, \epsilon)\leq k^*(y_0).
$$
(The latter inequality is valid due to the fact that $\delta > 0 $ can be arbitrary small). Thus (\ref{eq-est-from-above-18}) is established and the proof of the proposition is completed.
 $\ \Box$

\medskip

\medskip

{\bf Corollary 3.2.} {\it Let the assumptions of Proposition 3.1 be satisfied. Then the equality
\begin{equation}\label{e-SC-1}
v^*(y_0)=k^*(y_0)
\end{equation}
is valid if
\begin{equation}\label{e-SC-2}
d^*(y_0)=k^*(y_0)
\end{equation}
or if the equality (\ref{e-main-8-3}) is true.}

{\it Proof. } The proof follows from Proposition 2.3 and Proposition 3.1 or from Corollary 2.5 and Proposition 3.1 $\ \Box $

\medskip

In the next section, we will consider a class of systems for which  the equalities (\ref{e-SC-1}) and (\ref{e-SC-2}) are established to be valid.

\section{One special class of \lq\lq non-ergodic" systems}\label{Section-examples}

Let us introduce the following assumptions.

\medskip

ASSUMPTION III. There exists a continuously differentiable vector function
$F(y)=(F_i(y)), \ i=1,...,k,$ such that
\begin{equation}\label{e-CSO-3}
\nabla F_i(y)^Tf(y,u)= 0   ,\ i=1,...,k, \ \ \ \ \forall \ (y,u)\in \hat Y\times U
\end{equation}
where $\hat Y$ is a sufficiently large compact set.

\medskip

Define the set $Y_z$ by the equation
\begin{equation}\label{e-CSO-4}
Y_z:=\{y\in \R^m \ : \ F(y)=z\}.
\end{equation}
and define the set $Y$ as the union
\begin{equation}\label{e-SC-3}
Y=\bigcup_{z\in Z}Y_z,
\end{equation}
where $Z $ is some compact subset of $\R^k $. We assume that $Y $ is contained in $\hat Y $, that is, (\ref{e-CSO-3}) is satisfied for all $(y,u)\in Y\times U $. Therefore, each of the sets
 $Y_z, \ z\in Z$ and  the set $Y$ are forward invariant with respect to  system (\ref{e-CSO}).

  In addition to Assumption III, let us also introduce the following assumption.

  \medskip

ASSUMPTION IV. For any $z\in Z$ and for any
$y^1,y^2\in Y_z $, there exists a control $u(\cdot)$ ($u(t)\in U$) that steers  system
(\ref{e-CSO}) from $y^1$ to $y^2$ in finite time $ \mathcal{T}(y^1, y^2,z)\leq \mathcal{T}_0 $ ($\mathcal{T}_0$ being some positive constant).

\medskip

Due to Assumption III,  system (\ref{e-CSO}) is not ergodic on $Y$. However, Assumption IV  makes it  ergodic on each of $Y_z$ for $z\in Z$. To illustrate these assumptions, let us consider the following elementary example, in which they are readily verifiable.

\medskip

{\it Example.} Let $y(\tau)= (y_1(\tau), y_2(\tau))\in \R^2  $, \ $u(\tau)\in [-1,1]\in \R^1$, $\ f(y,u) = (f_1(y,u), f_2(y,u)) $ with
$\ f_1(y,u) = uy_2 $ and $\ f_2(y,u) = -uy_1 $. That is, system (\ref{e-CSO}) is of the form
$$
y_1'(t)= u(t)y_2(t), \ \ \ \ \ \ \  \ \ \ \ \ y_2'(t)= -u(t)y_1(t) .
$$
  It can be  seen that Assumption III
 is satisfied in this case with $F(y)=y_1^2 + y_2^2$ \ ($k=1$) and
$$
Y_z = \{(y_1,y_2) \ | \ y_1^2 + y_2^2 = z\} \ \  \forall \ z\geq 0.\
$$
Assuming that $Z= [a,b] $, where   $0<a<b$ are some constants, one can also see that,  with the use of the control $u(t)=1  $, any point in the set $Y_z$ (\ $z\in Z $)  can be reached from any other point of this set within a time interval that is less or equal than $2\pi $.
 Thus, Assumption IV is satisfied as well in this case.

 \medskip

The following proposition establishes  that the equalities (\ref{e-SC-1}) and (\ref{e-SC-2}) are valid for the class of systems satisfying Assumptions III and IV.\\

{\bf Proposition 4.1.} {\it Let Assumptions III and IV be satisfied. Then
\begin{equation}\label{e-CSO-further-6-0}
\lim_{T\rightarrow\infty}v_T(y_0)=k^*(y_0)= d^*(y_0)= \tilde{k}^*(z) \ \ \ \ \  \forall y_0\in Y_z,
\end{equation}
where $\tilde k^*(z)$  is the optimal value of the IDLP problem
\begin{equation}\label{e-CSO-further-6-1}
\tilde{k}^*(z)=  \min_{\gamma\in \mathcal{W}(z)}\int_{Y_z\times U}k(y,u)\gamma(dy,du),
\end{equation}
in which the minimization is over the set $\mathcal{W}(z)$ defined by the equation
\begin{equation}\label{limits-ergodic-W-Y-z}
 \mathcal{W}(z):= \left\{\gamma\in W \ : \ supp(\gamma)\in Y_z\times U \right\}
 \end{equation}
 ($W$ being defined in (\ref{limits-ergodic-W}) and $ supp(\cdot)$ standing for the support of the corresponding measure).
  Also,
 \begin{equation}\label{e-CSO-further-6-2}
cl W_1(y_0)= W_2(y_0)=  \mathcal{W}(z) \ \ \ \ \ \ \forall y_0\in Y_z,
\end{equation}
where $W_1(y_0) $ and $W_2(y_0) $ are defined in (\ref{limits-non-ergodic-2}) and (\ref{limits-non-ergodic-4})}.

\medskip



{\it Proof.} Note  that  we do not need to distinguish between the set $ \mathcal{W}(z)$ defined in (\ref{limits-ergodic-W-Y-z}) and the set defined by the equation
\begin{equation}\label{limits-ergodic-W-Y-z-1}
\left\{\gamma\in \mathcal{P}(Y_z\times U) \ : \ \int_{U\times Y_z}\nabla \phi(y)^Tf(u,y)\gamma(du,dy)=0 \ \ \forall \phi(\cdot)\in C^1 \right\},
\end{equation}
where $\ \mathcal{P}(Y_z\times U) $   stands for the space of probability measures defined on the Borel subsets of $Y_z\times U $. This set will also be denoted as $\mathcal{W}(z) $.

It can be established  that Assumptions III and IV imply that the following statement is valid (see, e.g., Proposition  3.3 in \cite{GR}):
{\it  For any $z\in Z$, given two arbitrary initial conditions  $y_0^1, y_0^2\in Y_z  $ and an arbitrary control $u^1(\cdot)\in \U$, there exists a control $u^2(\cdot)\in \U$ such that, for any continuous $q(u,y)$,
\begin{equation}\label{e-CSO-further-5}
|\frac{1}{T}\int_0^T q(y(t, y_0^1, u^1), u^1(t)dt - \frac{1}{T}\int_0^T q(y(t, y_0^2, u^2), u^2(t)dt |\leq \beta_q(T),
\end{equation}
for some $\beta_q(T) $ such that $\ \lim_{T\rightarrow\infty}\beta_q(T) = 0$.}

Due to the validity of this statement, from
Theorem 2.1(iii) and Proposition 4.1 in \cite{Gai8} it follows that
\begin{equation}\label{e-CSO-further-6}
\rho_H\left(\Gamma_T(y_0), \mathcal{W}(z)\right)\leq \beta(T)\ \ \ \ \forall y_0\in Y_z
\end{equation}
for some $\beta(T) $ such that $\ \lim_{T\rightarrow\infty}\beta(T) = 0$. By (\ref{e:occup-meas-def-eq-2}), the latter implies that
\begin{equation}\label{e-CSO-further-8}
\lim_{T\rightarrow\infty}v_T(y_0)=\tilde{k}^*(z)\ \ \ \ \forall y_0\in Y_z.
\end{equation}
Let us prove that
\begin{equation}\label{e-CSO-extra-2}
W_2(y_0)= \mathcal{W}(z) \ \ \ \ \ \ \forall y_0\in Y_z.
\end{equation}
To this end, let us first show that
\begin{equation}\label{e-CSO-extra-1}
W_2(y_0)\subset  \mathcal{W}(z) \ \ \ \ \ \ \forall y_0\in Y_z.
\end{equation}
Define the function
\begin{equation}\label{e-CSO-further-10}
  \Psi_z(y):= \sum_{i=1}^k (F_i(y)-z_i)^2.
\end{equation}
Note that, according to this definition,
\begin{equation}\label{e-SC-5}
\Psi_z(y)= 0\  \forall\ y\in Y_z , \ \ \ \ \ \ \Psi_z(y)> 0\  \forall\ y\in Y/Y_z
\end{equation}
and also
\begin{equation}\label{e-SC-4}
\nabla \Psi_z(y)^T f(y,u)= 2\sum_{i=1}^k (F_i(y)-z_i)\left(\nabla F_i(y)^T f(y,u)\right)= 0\ \ \ \ \forall \ (y,u)\in Y\times U.
\end{equation}
Take an arbitrary $\gamma\in W_2(y_0) $. Due to definition of $W_2(y_0) $ (see (\ref{limits-non-ergodic-4})), it implies that $\gamma\in W$ and that there exists a sequence
$\xi_l\in \mathcal{M}_{+}(Y\times U), \ l=1,2,..., $ such that
$$
\int_{Y\times U}(\Psi_z(y)-\Psi_z(y_0))\gamma(du,dy) \ =\ \lim_{l\rightarrow\infty}\int_{Y\times U}\nabla \Psi_z(y)^Tf(u,y)\xi_l(du,dy) \ = 0,
$$
where the equality to $0$ follows from (\ref{e-SC-4}). From this equality and from (\ref{e-SC-5}) it follows also that
$\ supp(\gamma)\subset Y_z $. Thus $\gamma\in \mathcal{W}(z) $ and the inclusion (\ref{e-CSO-extra-1}) is proved.

Take now an arbitrary $\gamma \in \mathcal{W}(z) $. By (\ref{e-CSO-further-6}), there exist
$T_l\rightarrow\infty $ and $\gamma_l\in  \Gamma_{T_l}(y_0)$ such that $\gamma_l\rightarrow\gamma $ as $l\rightarrow\infty$.
The fact that the measure $\gamma_l$ belongs to the set $\Gamma_{T_l}(y_0) $ means that this measure is generated by some control $u_l(\cdot)\in \U_{T_l}(y_0) $ and the corresponding solution $y_l(t)=y(t,y_0,u_l) $ of system (\ref{e-CSO}). Thus, the equality
(\ref{e-main-4}) is valid for any $\phi(\cdot)\in C^1 $. Proceeding now in exactly the same way as in the proof of Proposition 2.3, we obtain that $\gamma\in W_2(y_0) $. Consequently, $\ \mathcal{W}(z)\subset W_2(y_0) $ and by (\ref{e-CSO-extra-1}), the equality (\ref{e-CSO-extra-2}) is valid.

From (\ref{e-CSO-extra-2}) and from (\ref{limits-non-ergodic-3-1}), (\ref{e-CSO-further-6-1}) it follows that
\begin{equation}\label{e-CSO-further-7}
  d^*(y_0) = \tilde{k}^*(z) \ \ \ \ \forall y_0\in Y_z.
\end{equation}
To finalize the proof of (\ref{e-CSO-further-6-0}), we now only need to show
that
\begin{equation}\label{e-SC-4-1-5}
k^*(y_0) = \tilde{k}^*(z) \ \ \ \ \forall y_0\in Y_z.
\end{equation}
From (\ref{limits-non-ergodic-dual-4}) and (\ref{e-CSO-further-7})
 it follows that $k^*(y_0) \geq\tilde{k}^*(z) $. On the other hand, by Lemma 4.2 (see below), the equality (\ref{e-main-8-3}) is true. Therefore, by Corollary 2.5 (and
thanks to (\ref{e-CSO-further-8})), $\tilde{k}^*(z)\geq k^*(y_0)  $. Thus, (\ref{e-SC-4-1-5}) is true and (\ref{e-CSO-further-6-0}) is proved.

The fact that
$cl W_1(y_0) = \mathcal{W}(z) \ \forall \ y_0\in Y_z $ follows from the fact that (\ref{e-SC-4-1-5}) is valid with the use of   any $k(y,u)$ in (\ref{limits-non-ergodic-1})
and (\ref{e-CSO-further-6-1}), and also from the fact that the sets $W_1(y_0) $  and $\mathcal{W}(z) $ are convex.
Since  (\ref{e-CSO-extra-2}) has been already established, the proof of the proposition is completed.
 $\ \Box$

{\bf Lemma 4.2.} {\it If Assumptions III and IV are satisfied, then the equality (\ref{e-main-8-3}) is true.}

{\it Proof.} By (\ref{e:occup-meas-def-eq-per-1}), (\ref{e-main-8-2}) and (\ref{e-CSO-further-8}), $\  v_{per}(y_0)\geq  \tilde{k}^*(z) $.
Thus,
 to prove (\ref{e-main-8-3}), it is sufficient to prove that
\begin{equation}\label{e-CSO-further-17}
 \tilde{k}^*(z)\geq v_{per}(y_0)\ \ \ \forall y_0\in Y_z.
\end{equation}
 Let $x_0\in Y_z$. By (\ref{e-CSO-further-8}), for any sequence $T_l, l=1,2,...$, $\ T_l\rightarrow\infty $, there exist $\ u^l(\cdot)\in \U_{T_l} $, $\ y^l(t)= y(t,y_0,u^l)  $  such that
\begin{equation}\label{e-CSO-further-13-1}
\lim_{\l\rightarrow\infty}\frac{1}{T_l}\int_0^{T_l}k(y^l(t), u^l(t))dt = \tilde{k}^*(z).
\end{equation}
By Assumption IV, there exists a control $\ u(t)\in U$ defined on an interval $ \ t\in [T_l, T_l+\Delta_l]$ (with $\ 0\leq \Delta_l\leq \mathcal{T}_0 $) such that, with the use of this control, system (\ref{e-CSO}) will be steered from the point $y^l(T_l) $ at $t=T_l$ to the point $y_0$ at $t=T_l+\Delta_l$. Denote by $\tilde{u}^l(\cdot) $ the control that is equal to $u^l(\cdot)$ on the interval $[0, T_l) $ and equal to the \lq\lq steering control" on the interval $[T_l, T_l + \Delta_l] $. Denote also by $\tilde{y}^l(\cdot) $ the corresponding solution of (\ref{e-CSO}). The definition of the pair $(\tilde{y}^l(\cdot),\tilde{u}^l(\cdot)) $ can be extended to the infinite time horizon as $\mathcal{T}_l$-periodic pair with $\mathcal{T}_l= T_l+\Delta_l$. Therefore, by (\ref{e-main-8-1}),
\begin{equation}\label{e-CSO-further-13-2}
  \frac{1}{\mathcal{T}_l}\int_0^{\mathcal{T}_l}k(\tilde y^l(t), \tilde u^l(t))dt\geq v_{per}(y_0) , \ \ l=1,2,...
\end{equation}
On the other hand,
$$
\left|\frac{1}{T_l}\int_0^{T_l}k(y^l(t), u^l(t))dt - \frac{1}{\mathcal{T}_l}\int_0^{\mathcal{T}_l}k(\tilde y^l(t), \tilde u^l(t))dt\right|
$$
$$
\leq \left|\frac{1}{T_l}\int_0^{T_l}k(y^l(t), u^l(t))dt - \frac{1}{\mathcal{T}_l}\int_0^{T_l}k( y^l(t),  u^l(t))dt\right|
$$
$$
+ \frac{1}{\mathcal{T}_l}\int_{T_l}^{T_l+\Delta_l}|k(\tilde y^l(t), \tilde u^l(t))|dt\leq \frac{M_k\mathcal{T}_0}{\mathcal{T}_l},
\ \ \ {\rm where} \ \ \ M_k:=\max_{(y,u)\in Y\times U}|k(y,u)|.
$$
The latter and (\ref{e-CSO-further-13-1}) imply
$$
\lim_{l\rightarrow\infty}\frac{1}{\mathcal{T}_l}\int_0^{\mathcal{T}_l}k(\tilde y^l(t), \tilde u^l(t))dt = \tilde{k}^*(z),
$$
which, in turn, implies (\ref{e-CSO-further-17}) (due to (\ref{e-CSO-further-13-2})). \ $\ \Box $

\section{Proofs of some duality results}\label{Section-Duality-proofs}

 Let $(C^1)^* $ stand for the space of continuous functionals on the space of smooth functions $C^1$, the latter being considered as the normed vector space with norm $||\phi(\cdot)||:= \max_{y\in Y}|\phi(y)| + \max_{y\in Y}||\nabla \phi(y)|| $ for any $\phi(\cdot)\in C^1 $.
  Define a linear operator $\A(\cdot): \M(Y\times U) \times \M(Y\times U)\mapsto \R^1\times (C^1)^* \times (C^1)^* $  as follows: for any $(\gamma, \xi)\in \M(Y\times U)\times \M(Y\times U) $
\begin{equation}\label{e-Duality-1}
\A(\gamma, \xi):= \left(\int_{Y\times U}\gamma(dy,du), \ a_{(\gamma, \xi)}, \ b_{\gamma}\right),
\end{equation}
where $a_{(\gamma, \xi)}, \ b_{\gamma}\in (C^1)^*$ are defined by the equation: $\ \forall \  \phi(\cdot)\in C^1$,
$$
 \ a_{(\gamma, \xi)}(\phi) := -\left\{\int_{Y\times U} (\phi(y_0)-\phi(y))\gamma(dy,du)
 + \int_{Y\times U}\nabla\phi (y)^Tf(y,u)\xi(dy,du)\right\},
 $$
 \vspace{-0.4cm}
\begin{equation}\label{e-Duality-2-1}
\  b_{\gamma}(\phi):= - \left\{\int_{Y\times U}\nabla\phi (y)^Tf(y,u)\gamma(dy,du)\right\}.
\end{equation}
In this notation, the set $\Omega(y_0)$ (defined in (\ref{non-ergodic-Omega})) can be rewritten as follows
\begin{equation}\label{e-Duality-2-0-11}
 \Omega(y_0)= \{(\gamma , \xi)\in \M_+(Y\times U)\times \M_+(Y\times U)\ : \ \A(\gamma, \xi)= (1, {\bf 0}, {\bf 0})\}
\end{equation}
where  ${\bf 0} $ stands for the zero element of $(C^1)^*$. Also,
problem (\ref{limits-non-ergodic}) takes the form
\begin{equation}\label{e-Duality-2-0}
\inf_{(\gamma , \xi)\in \Omega(y_0)}\langle k, \gamma  \rangle\ := k^*(y_0)
\end{equation}
where
 $\langle \cdot, \gamma  \rangle $ (also,
$\langle \cdot, \xi  \rangle $ in the sequel) denoting the integral of the corresponding function over $\gamma$ (respectively, over $\xi$).

Note that for any $(\mu, \psi(\cdot), \eta(\cdot))\in \R^1\times C^1\times C^1 $,
$$
\langle A(\gamma, \xi), (\mu, \psi, \eta)\rangle =
\mu \int_{Y\times U}\gamma(dy,du) + a_{(\gamma, \xi)}(\psi)
+ b_{\gamma}(\eta)
$$
$$
=  \int_{Y\times U}\left(\mu - (\psi(y_0)-\psi(y)) - \nabla \eta(y)^Tf(y,u)\right)\gamma(dy,du)
$$
\begin{equation}\label{e-Duality-2-2-0}
- \int_{Y\times U} \nabla \psi(y)^Tf(y,u)\xi(dy,du).
\end{equation}

Define now the linear operator $\A^*(\cdot): \R^1\times C^1\times C^1\mapsto C(Y\times U)\times C(Y\times U)\subset \M^*(Y\times U)\times  \M^*(Y\times U) $ in such a way that, for any $(\mu, \psi(\cdot),\eta(\cdot))\in \R^1\times C^1\times C^1 $,
\begin{equation}\label{e-Duality-3-1}
\A^*(\mu, \psi, \eta)(y,u):= \left(\mu - (\psi(y_0)-\psi(y)) - \nabla \eta(y)^Tf(y,u), \ \nabla \psi(y)^Tf(y,u)\right).
\end{equation}
Thus
$$
\langle (\gamma, \xi), \A^*(\mu, \psi, \eta) \rangle =  \int_{Y\times U}\left(\mu - (\psi(y_0)-\psi(y)) - \nabla \eta(y)^Tf(y,u)\right)\gamma(dy,du)
$$
$$
 - \int_{Y\times U} \nabla \psi(y)^Tf(y,u)\xi(dy,du) = \langle A(\gamma, \xi), (\mu, \psi, \eta)\rangle  .
$$
That is, the operator $ \A^*(\cdot) $ is the adjoint of  $ \A(\cdot) $. The problem dual to (\ref{e-Duality-2-0})
is of the form (see \cite{And-1} and \cite{And-2})
\begin{equation}\label{e-Duality-4-0}
\sup_{(\mu, \psi(\cdot), \eta(\cdot))\in \R^1\times C^1\times C^1} \mu = d^*(y_0)
\end{equation}
$$
\ s.\ t.  \ \ \ \ \ \ \ \ \ \ \   \ \ \ \ \ \ \ \ \ \ \ \ \ \ \ \ \ \ \ \   \ \ \ \ \ \ \ \ \
$$
\begin{equation}\label{e-Duality-4-0-1}
 - \A^*(\mu, \lambda, \nu)(y,u) + (k(y,u), 0)\geq (0,0) \ \ \forall (y,u)\in Y\times U,
\end{equation}
the latter being equivalent to (\ref{limits-non-ergodic-dual}).

\bigskip

{\it Proof of Lemma  2.2.} Let
$$
H:= \Big\{\left( \A(\gamma, \xi), \int_{Y\times U}k(y,u)\gamma(dy,du) + r\right)\ :
$$
\vspace{-0.6cm}
\begin{equation}\label{limits-non-ergodic-pert-dual-4}
\ (\gamma, \xi)\in \M_+(Y\times U)\times \M_+(Y\times U), \ r\geq 0\}\subset \R^1\times (C^1)^*\times (C^1)^*\times \R^1\ \Big\} ,
\end{equation}
and let $\bar{H}$ stand for the closure of $H$ in the weak$^*$ topology of $\R^1\times (C^1)^*\times (C^1)^*\times \R^1$.
Consider the problem
\begin{equation}\label{limits-non-ergodic-pert-dual-4-1}
\inf \{\theta \ | \ (1, {\bf 0}, {\bf 0}, \theta)\in \bar{H} \}:= k_{sub}^*(y_0).
\end{equation}
Its optimal value $k_{sub}^*(y_0)$ is called the subvalue of the IDLP problem (\ref{e-Duality-2-0}). Let us show that the optimal value of (\ref{limits-non-ergodic-3}) is equal to the subvalue. In fact,
as can be readily seen, $\left(1, {\bf 0}, {\bf 0},  \int_{Y\times U}k(y,u)\gamma(dy,du)\right)\in \bar{H} $  if $\gamma\in W_2(y_0)$. Consequently,
\begin{equation}\label{limits-non-ergodic-pert-dual-5}
k_{sub}^*(y_0)\leq \min_{\gamma\in W_2(y_0)}\int_{Y\times U}k(y,u)\gamma(dy,du).
\end{equation}
From the fact that $k_{sub}^*(y_0)$ is defined as  the optimal value in (\ref{limits-non-ergodic-pert-dual-4-1}) it follows that there exists  a sequence $(\gamma_l,\xi_l)\in \M_+(Y\times U)\times \M_+(Y\times U) $  such that $\A(\gamma_l,\xi_l) $ converges (in weak$^*$ topology)
to $(1, {\bf 0}, {\bf 0})$, with   $\int_{Y\times U}k(y,u)\gamma_l(dy,du) $ converging to $k_{sub}^*(y_0)$ as $l$ tends to infinity. That is (see (\ref{e-Duality-1})),
$$
\int_{Y\times U}\gamma_l(dy,du)\rightarrow 1, \ \ a_{(\gamma_l, \xi_l)}\rightarrow {\bf 0},\ \  b_{\gamma_l}\rightarrow {\bf 0},
$$
 \vspace{-0.4cm}
$$
\ \ \int_{Y\times U}k(y,u)\gamma_l(dy,du)\rightarrow k_{sub}^*(y_0).
$$
Without loss of generality, one may assume
that $\gamma_l$ converges in weak$^*$ topology
 to a  measure $\gamma$ that satisfies the relationships
 $$
 \int_{Y\times U}\gamma(dy,du)=1,  \ \  b_{\gamma}= {\bf 0}\ \ \ \ \Rightarrow \ \ \ \  \gamma\in W.
 $$
 Also,
 $
 \  a_{(\gamma, \xi_l)}\rightarrow {\bf 0}
 $ and $\int_{Y\times U}k(y,u)\gamma(dy,du)= k_{sub}^*(y_0) $. That is, $\gamma\in W_2(y_0)$ and therefore,
 \begin{equation}\label{limits-non-ergodic-pert-dual-5-1}
\min_{\gamma\in W_2(y_0)}\int_{Y\times U}k(y,u)\gamma(dy,du)\leq k_{sub}^*(y_0).
\end{equation}
 Thus the optimal value of (\ref{limits-non-ergodic-3}) is equal to the subvalue.
 To complete the proof, it is sufficient to note that the subvalue of an IDLP problem is equal to the optimal value of its dual provided that the former is bounded (see, e.g., Theorem 3 in \cite{And-1}). That is, $k_{sub}^*(y_0) = d^*(y_0). $ $\ \Box$

\bigskip

In the notation of this section, problem (\ref{limits-non-ergodic-pert-eps}) takes the form
\begin{equation}\label{e-Duality-3}
\inf_{(\gamma,\xi) \in \Omega(y_0)}\langle (k+2\epsilon, M\epsilon ), (\gamma , \xi) \rangle\ = k^*(y_0,\epsilon),
\end{equation}
while its dual (\ref{limits-non-ergodic-dual-eps-d}) takes the form
\begin{equation}\label{e-Duality-4-0-2}
\sup_{\mu, \psi(\cdot), \eta(\cdot)\in \R^1\times C^1\times C^1} \mu = d(y_0,\epsilon)
\end{equation}
$$
\ s.\ t.  \ \ \ \ \ \ \ \ \ \ \   \ \ \ \ \ \ \ \ \ \ \ \ \ \ \ \ \ \ \ \   \ \ \ \ \ \ \ \ \
$$
\begin{equation}\label{e-Duality-4-0-2}
 - \A^*(\mu, \lambda, \nu)(y,u) + (k(y,u)+2\epsilon, M\epsilon)\geq (0,0) \ \ \forall (y,u)\in Y\times U .
\end{equation}

\bigskip

{\bf Lemma 5.1.} {\it Let $\Omega(y_0)\neq\emptyset $. Then the optimal values of (\ref{e-Duality-3}) and (\ref{e-Duality-3})-(\ref{e-Duality-4-0-2}) are equal for any $\epsilon > 0 $. That is, (\ref{eq-est-from-above-16}) is valid. }

\bigskip

{\it Proof.} By Theorem 6 in \cite{And-1}, to prove the validity of (\ref{eq-est-from-above-16}), it is sufficient to establish
that the set $D$,
\begin{equation}\label{Lemma-Anderson-1}
D:=\{ \left(\mathcal{A}(\gamma , \xi),\ \langle (k+2\epsilon, M\epsilon ), (\gamma , \xi) \rangle\ \right) \ : \ (\gamma , \xi)\in
{\cal M}_+(Y\times U)\times {\cal M}_+(Y\times U)\}
\end{equation}
is closed in weak$^*$ topology of $ \R^1\times (C^1)^*\times (C^1)^*\times \R^1$. The proof of this is
similar to the proof of Theorem 12 in \cite{And-1}. It is based on the following two properties of the problem.

{\it Property 1.} The set  ${\cal M}_+(Y\times U))\times {\cal M}_+(Y\times U) $ has a compact base. That is (see \cite{And-1}),
\begin{equation}\label{Lemma-Anderson-2}
{\cal M}_+(Y\times U)\times {\cal M}_+(Y\times U) = \{\lambda (\gamma , \xi) : \ (\gamma , \xi)\in \mathcal{L}, \ \lambda\geq 0  \},
\end{equation}
where
\begin{equation}\label{Lemma-Anderson-2}
\mathcal{L}:= \{(\gamma , \xi)\in  {\cal M}_+(Y\times U)\times {\cal M}_+(Y\times U) \ :
\end{equation}
$$
\int_{Y\times U} \gamma(dy,du)  + \int_{Y\times U} \xi(dy,du) = 1 \},
$$
with $\mathcal{L}$ being a weak$^*$ compact subset of ${\cal M}(Y\times U)\times {\cal M}(Y\times U)$.

\medskip

{\it Property 2.} For $(\gamma, \xi)\in  {\cal M}_+(Y\times U)\times {\cal M}_+(Y\times U)$, the equalities
\begin{equation}\label{Lemma-Anderson-3}
\mathcal{A}(\gamma , \xi) = (0, {\bf 0}, {\bf 0}), \ \ \ \ \langle (k+2\epsilon, M\epsilon ), (\gamma , \xi) \rangle = 0
\end{equation}
can be valid only if $\gamma = 0 $ and $\xi = 0$, this being readily verifiable due to the fact that a part of the relationships  (\ref{Lemma-Anderson-3}) are the following equalities:
$$
\int_{Y\times U} \gamma(dy,du) = 0, \ \ \ \ \ \ \int_{Y\times U}(k(y,u)+2\epsilon) \gamma(dy,du)  + M\epsilon \int_{Y\times U} \xi(dy,du) = 0.
$$
Let us now prove that $D$ is closed. Let $(\gamma_l , \xi_l)\in {\cal M}_+(Y\times U)\times {\cal M}_+(Y\times U) $ be such that
\begin{equation}\label{Lemma-Anderson-4}
\mathcal{A}(\gamma_l , \xi_l) - {\bf z}\rightarrow (0, {\bf 0}, {\bf 0}), \ \ \ \ \langle (k+2\epsilon, M\epsilon ), (\gamma_l , \xi_l) \rangle\ - \beta\rightarrow 0 \ \ \ {\rm as} \ \ \ l\rightarrow\infty,
\end{equation}
where
${\bf z}\in  \R^1\times (C^1)^*\times (C^1)^*  $ and $\beta\in \R^1$. By (\ref{Lemma-Anderson-2}), $(\gamma_l , \xi_l) = \lambda_l (\bar \gamma_l , \bar \xi_l) $, where $\gamma_l\geq 0 $ and $\ (\bar \gamma_l , \bar \xi_l)\in \mathcal{L} $. Due to compactness of $\mathcal{L} $, one may assume (without loss of generality) that  $\ (\bar \gamma_l , \bar \xi_l)\rightarrow (\bar \gamma , \bar \xi)
\in \mathcal{L} $. Note that, due to Property 2, the sequence $\lambda_l$ is bounded. Indeed, assuming that this is not the case and there exists a subsequence $\{l'\} $ of $\{l\} $ such that $\lambda_{l'}\rightarrow\infty $ as $l'\rightarrow\infty $ , one would obtain (via substitution of $\ \lambda_{l'} (\bar \gamma_{l'} , \bar \xi_{l'}) $ into (\ref{Lemma-Anderson-4}) and passing to the limit with $l'\rightarrow\infty $) that
$$
\ \mathcal{A}(\bar \gamma_{l'} ,\bar \xi_{l'})- \frac{1}{\lambda_{l'}}{\bf z}\rightarrow (0, {\bf 0}, {\bf 0}), \ \ \
  \langle (k+2\epsilon, M\epsilon ), (\bar \gamma_{l'} , \bar \xi_{l'}) \rangle\ -
\frac{1}{\lambda_{l'}}\beta\ \rightarrow
0 \ \ \ {\rm as} \ \ \ l\rightarrow\infty.
$$
$$
\Rightarrow \ \ \ \ \ \mathcal{A}(\bar\gamma , \bar \xi) = (0, {\bf 0}, {\bf 0}), \ \ \ \ \langle (k+2\epsilon, M\epsilon ), (\bar\gamma , \bar \xi) \rangle = 0.
$$
According to Property 2, the latter implies that $\bar\gamma = 0 $ and $\bar\xi = 0$. This contradicts to the fact that $\ (\bar \gamma , \bar \xi)\in \mathcal{L} $. Thus the sequence $\{\lambda_l\} $ is bounded and therefore one may assume (without loss of generality) that $\lambda_l\rightarrow \lambda  $ as $l\rightarrow\infty $. Consequently,  $(\gamma_l, \xi_l)\rightarrow \lambda (\bar\gamma , \bar \xi) $
as $l\rightarrow\infty $. Denoting $(\lambda\bar\gamma , \lambda\bar \xi):=(\gamma , \xi) $,  one obtains (by (\ref{Lemma-Anderson-4}))
$$
\mathcal{A}(\gamma , \xi) = {\bf z}, \ \ \ \langle (k+2\epsilon, M\epsilon ), (\gamma , \xi) = \beta \ \ \ \ \Rightarrow \ \ \ \
({\bf z},\beta)\in D.
$$
This proves that $D$ is closed. $\ \Box$

\end{document}